\documentclass[12pt,fleqn,a4paper]{article}
 \usepackage{amssymb}
 \usepackage{amsfonts}
 \usepackage{amsmath}
\usepackage{graphicx}
\usepackage{latexsym}
\usepackage{enumerate}
 \usepackage{dsfont}
\usepackage{stmaryrd}
 \usepackage{mathrsfs}
 \usepackage{multicol}

 \newtheorem{theorem}{\sc Theorem}[section]
 
 \newtheorem{lemma}[theorem]{\sc Lemma}
 
 \newtheorem{corollary}[theorem]{\sc Corollary}

 \numberwithin{equation}{section}
 \DeclareMathAlphabet{\mathsfn}{OT1}{cmss}{x}{n}
\DeclareMathAlphabet{\mathsfb}{OT1}{cmss}{bx}{n}

  \newcommand{\bfm}{\boldsymbol}
  \newcommand{\alg}{\bfm}

 \newcommand{\bg}{\begin}
 \newcommand{\e}{\end}
 \newcommand{\mbb}{\mathsf}
 \newcommand{\be}{\bg{enumerate}}
 \newcommand{\ee}{\e{enumerate}}
 \newcommand{\bi}{\bg{itemize}}
 \newcommand{\ei}{\e{itemize}}

 \newcommand{\Lra}{\Leftrightarrow}
 
 \newcommand{\beqn}{\bg{eqnarray*}}

 \newcommand{\eeqn}{\e{eqnarray*}}
 \newcommand{\smsm}{\smallsetminus}
 
 \newcommand{\no}{\noindent}
 
 \newcommand{\A}{\alg{A}}
 
 \newcommand{\LLn}{L^{\omega}_{n+1}}

 \newcommand{\Ap}{A^{p+1}_{n+1}}
 \newcommand{\aAp}{\bfm{A^{p+1}_{n+1}}}

 \newfont{\gros}{cmr10 at 13pt}

 \newcommand{\tdiv}{\stackrel{\odot\ }{\to}}

\def\le{\leqslant}
\def\ge{\geqslant}

\def\brl{bounded resi\-duated lattice}

\def\BRL{\mbb{BRL}}
\def\WL{\mbb{WL}}
\def\WLn{\mbb{WL_\mbb{n+1}}}

\def\Em{\mbb{E_{m}}}

\def\EMm{\mbb{EM}_{\mbb{m}}}

\def\EMp{\mbb{EM}_{\mbb{p}}}

\begin{document}

\title{Families of  local involutive integral  residuated lattice-ordered commutative monoids admitting  Boolean  term.}

\author{by \sc Antoni Torrens}
\date{}
\maketitle

\begin{abstract}
  We  present a family of local  involutive integral bounded residuated lattice-ordered commutative monoids (involutive  residuated lattices, for short) having  Boolean term,  radical term (see \cite{CT12} and \cite{T23}), and satisfying GAP (Generalized Appel property) (see \cite{T23}).
 The construction of this family is based in the examples given in \cite[Subsection 5.1]{CT12}, by taking $\bfm{\LLn}$ in place of $\bfm{L_{n+1}}$ and $\bfm{L_{p+1}}$ in place of the two element Boolean algebra.
 In some sense this paper is a companion of the paper \cite{T23}, in which the General Apple Property (GAP) and Boolean terms are studied.
\end{abstract}

\section{Purpose and  preliminaries }

 In this paper we give the family $\{\aAp=\langle\Ap,\odot,\tdiv\,,\land,\lor,\bot,\top \rangle :n,p>0\}$ of algebras of type $(2,2,2,2,0,0)$ such that for any $n,p>0$ the algebra
 $\aAp$
 is an \emph{involutive integral residuated lattice-ordered commutative monoid} (\emph{involutive residuated lattice for short}), that is,
 \bi
 \item $\aAp$ is and {\em integral bounded residuated lattice-ordered commutative monoid (bounded residuated lattice, for abbreviate)}, i.e.,
     \bi
 \item $\langle \Ap; \odot , \top\rangle$ is a commutative monoid,
 \item $\langle \Ap; \lor, \land,\bot, \top\rangle$ is a bounded
 lattice with smallest element $\bot$, and greatest element $\top,$
 \item  the following \emph{residuation
 condition} is satisfied:  for any $a,b,c\in \Ap$
 \bg{align}\label{Eq:resd-cond}
 a \odot  b \le  c &
 \mbox{ if and only if } a \le b\tdiv c,\end{align}
 where  $\le $ is the order given by the
 lattice structure.
 \ei
 \item And $\aAp$  is \emph{involutive}:

 for all $a\in \Ap$, $a=(a\tdiv\bot)\tdiv\bot$ (or  $\sim \sim\!a=a$, where $\sim x=x\tdiv\bot$).\smallskip
 \ei

 For any $n,p>0$, $\aAp$ is \emph{local}, that is, it has a unique maximal proper congruence relation (a unique maximal proper implicative filter), i.e., $\langle Con(\aAp),\subseteq \rangle$ has  only one coatom; it is also \emph{directly indecomposable}, or equivalently, by [KO01, Proposition 1.5], it has only two Boolean (complemented) elements, namely $\bot$ and $\top$.
 Moreover,   $\aAp$  admits a Boolean term satisfying condition $(3.1)$ of Theorem 3.10 of \cite{T23}, and so by  \cite[Corollary 3.14]{T23}, its generated variety $HSP(\aAp)$ has GAP, that is,  directly indecomposable and local  coincide in $HSP(\aAp)$. We also see that this variety admits  radical term.\smallskip

  To describe the Boolean  terms mentioned above we must recall some definitions.
  In the algebraic  language $\langle \ast,\to,\land,\lor,\bot\top\rangle$ of type $(2,2,2,2,0,0)$ we define the following terms:
 \bi
 \item $\neg x= x\to\bot$,
 \item $x\oplus y =\neg (\neg x\ast \neg x)$,
 \item $x^0=\top$ and $0.x=\bot$,
 \item for any  integer $n>0$,  $x^{n+1}=x\ast x^n$ and $(n+1).x=x\oplus n.x$.
 \ei
 These terms are used in the literature   to define
 several subvarieties of the variety $\BRL$ of bounded residuated lattices. For example, if $m>0$,   then  the varieties  $\Em$, $\EMm$ and $\WL_\mbb{m}$  are the subvarieties of $\BRL$  given respectively  by the equations  $x^{m+1}\approx x^m$, $x\lor \neg x^m\approx\top$ and $m.x\lor m.\neg x\approx \top$ (see for example \cite{GJKO}, \cite{K04} \cite{CT12},\cite{T16} and \cite{T23}).
 Arithmetic properties of the interpretations of these terms in bounded residual lattices can be found in these papers and in the references given therein. In our case for any $m,p>0$   the interpretation of $\ast$ over $\aAp$ is $\odot$ (i.e. $\odot=\ast^{\aAp}$) and the interpretation of $\to$ over $\aAp$ is $\tdiv$ (i.e. $\tdiv\,=\to^{\aAp}$). \smallskip

 In Theorem \ref{T:Bterm} we show that $(n+1).x^{\max\{n+1,p\}}$ is  Boolean term and radical term for $\aAp$, and so it is Boolean term for the variety $HSP(\aAp)$.  This term satisfies the condition $(3.1)$ of \cite[Theorem 3.10]{T23}, hence $HSP(\aAp)$ satisfies GAP, and so its  directly indecomposable members are local.  We also see that the class $HSP(\aAp)_{\mbb{SS}}$ of semisimple members of  $HSP(\aAp)$ form a variety, actually  we show that $HSP(\aAp)_{\mbb{SS}}\subseteq \EMp$, then, by \cite[Theorem 3.17]{T23}, $HSP(\aAp)$ admits a boolean term that is also radical term.
 \smallskip

 For any $n,m>0$,  $\aAp$ is constructed  from the MV-algebras (or bounded Wajsberg hoops) $\bfm{\LLn}$ and $\bfm{L_{p+1}}$. To fix notation we are going to recall  how this MV-algebras can be obtained (see for example \cite{CDM00},\cite{GMT99} and \cite{GM05}).
 We consider  MV-algebras as \brl s, and we assume that the reader is familiar with these structures.\smallskip

 Let $\langle\mathds{Z};+,0;\le\rangle$ be the totally ordered group of integers, and let $\mathds{Z}\stackrel{\to}{\otimes}\mathds{Z}=\langle \mathds{Z}\times\mathds{Z},+,(0,0);\preccurlyeq  \rangle$ be  the lexicographic product of $\mathds{Z}$ by $\mathds{Z}$ whose order $\preccurlyeq$ is:
\[(m,r)\preccurlyeq (k,s) \mbox{ iff }\left\{\begin{array}{ll}
 m<k, & \mbox{ or   } \\
  m=k & \mbox{ and }r\le s\end{array}\right.\]

\no For any  integer $n>0$ we consider the MV-chain $\bfm{\LLn}=\bfm{\Gamma}(\mathds{Z}\stackrel{\to}{\otimes}\mathds{Z},(n,0))$, obtained from the totally ordered grup $ \mathds{Z}\stackrel{\to}{\otimes}\mathds{Z}$ with the strong unit $(n,0)$; considered as  a Wajsberg algebra (bounded Wajsberg Hoop), i.e.,
 \[\bfm{\LLn}=\left\langle L^{\omega}_{n+1}=[(0,0),(n,0)];\ast,\to,\land,\lor,(0,0),(n,0)\right\rangle;\]
where
 \bi
 \item
 $\LLn=\{(m,r)\in \mathds{Z}\times\mathds{Z}:(0,0)\preccurlyeq(m,r)\preccurlyeq(n,0) \},$
  \item its lattice order $\le_{\bfm{\LLn}}= \ \preccurlyeq\! |_{\LLn}$ is total,
  \item $(m,r)\ast (k,s)=\max\{(0,0),[(m,r)+ (k,s)]-(n,0)\}$\footnote{Here and  in what follows ``$\max$" and ``$\min$" are relative  to order $\preccurlyeq$ of $\mathds{Z}\stackrel{\to}{\otimes}\mathds{Z}$},
      \item $(m,r)\to (k,s)=\min\{(n,0),[(n,0)-(m,r)]+(k,s)\}$.
      \ei
  Moreover,
  \bi
  \item $(m,r)\ast(k,s)\le_{\bfm{\LLn}}(l,t)$ iff $(k,s)\le_{\bfm{\LLn}}(m,r)\to(l,t)$, and
  \item it is \emph{involutive}, indeed, if for any $(m,r)\in\LLn$ $\neg (m,r)=: (m,r)\to (0,0)$, then for any $(k,s)\in\LLn$,  $(m,r)\le_{\bfm{\LLn}}(k,s)$ iff $\neg (k,s)\le_{\bfm{\LLn}}\!\neg (m,r)$.  And we have $ (m,r)\to (k,s)= \neg\,[(m,r)\ast\neg(k,s)]$.
  \ei
   We also consider $L_{n+1}=\{\bfm{m}=(m,0)\in \mathds{Z}\times \mathds{Z}: 0\le m\le n\}$ which is universe of $\bfm{L_{n+1}}$  subalgebra of $\bfm{\LLn}$.  $\bfm{L_{n+1}}$   is also a copy of the  totally ordered Wajsberg hoop
  with $n+1$ elements with lower bound (see \cite{BF00}). In par\-ticular, $\bfm{L_{2}}$ is a copy of the two element Boolean algebra, in which  $\ast=\land$ and $\neg x$ is the complement of $x$.\smallskip

  Given $n,p>0$ to define $\aAp$ we procedure as follows:
  \be[(1)]
  \item As universe  $\Ap$ we take a suitable subset of $L^{\omega}_{n+1}\times L_{p+1}$.
  \item We define the bounded lattice reduct  $\langle \Ap,\land,\lor,\bot,\top\rangle$ by means of the lattice partial  order $\le$.
 \item With the help of  the operations $\ast$ and $\to$ of $\bfm{\LLn}$ and $\bfm{L_{p+1}}$, we define the monoidal operation $\odot$ on $\Ap$ and we show that $\langle \Ap;\odot,\top; \le\rangle$ is an integral  partial-ordered commutative monoid.
  \item We define an unary operation $\sim$, which is an involution ($x\le y$ if and only if $\sim y\le \sim x$), and we show that it   satisfies $x\odot y=\bot$ if and only if $x\le\sim y$.
  \item Finally, we define the operation $x\tdiv y=:\sim(x\odot \sim y)$, and   by Lemma \ref{L:Monid-invo} given  bellow it satisfies the residuation condition.
  \ee

  For item (5) we need a preliminary result that can be deduced from the results given in \cite{HR07}. However, we include an easy proof.
 \bg{lemma}\label{L:Monid-invo} Let $\langle M; \ast,\top;\le,\rangle$ be a commutative integral partial ordered monoid, and let $\neg:M\to M$ be an involution, i.e., for all $a,b\in M$ $a\le b$ iff $\neg b\le \neg a$.   If  for any $a,b\in M$  $a\ast b=\neg\top\mbox{ iff } a\le \neg b$,  then for any $a,b,c\in M$
 $a\ast b\le c \mbox{ iff }  b\le \neg(a\ast \neg c)$
 \e{lemma}
 {\em Proof:}
 Let $a,b$ and $c$ be arbitrary elements of $M$. Since $\neg$ is an involution then $\neg\neg a=a$. Therefore, taking account that $\ast$ is associative and commutative in $M$, we have
 \[
 a\ast b\le c \Lra  (a\ast b)\ast \neg c\le \neg\top \Lra b\ast (a\ast \neg c) \le \neg\top
 \Lra b\le \neg(a\ast \neg c).\hspace{0,6cm}\Box
 \]

  We assume that the reader is familiar wiht bounded integral commutative residuated lattices, MV-algebras and Wajsberg hoops, see  for example \cite{CDM00} \cite{BF00}, \cite{FRT84}, \cite{GJKO}  and the references given therein.
 The results on universal algebras used in this paper can be found in the books \cite{BS81} and \cite{Be12}.

\section{Building de family $\{\aAp :n,p>0\}$}

 For any pair of integers $n,p>0$, we consider the algebra
  \[\aAp=\left\langle \Ap ;\odot,\tdiv,\land,\lor,\langle (n,0),\bfm{0}\rangle,\langle (n,0),\bfm{p}\rangle\right\rangle\]
   of type $\langle 2,2,2,2,0,0\rangle$ defined as follows
 \subsection*{The universe $\Ap$}
 Its universe $\Ap$ is the following subset of $\LLn\times L_{p+1}$:  \[A^{p+1}_{n+1}=\big(\LLn\times\{\bfm{0},\bfm{p}\}\big)\cup \big(L^{\omega}_{n}\times (L_{p+1}\smsm\{\bfm{0},\bfm{p}\})\big).\]
In other words\smallskip

 \no $A^{p+1}_{n+1}=\big(\LLn\times L_{p+1}\big)\!\smsm\!\big(\{(m,r):(n-1,0)\prec(m,r)\preccurlyeq(n,0) \}\times \big(L_{p+1}\!\smsm\!\{\bfm{0},\bfm{p}\})\big)$.\medskip

 \no Then an   element of $\Ap$  is an  ordered pair $\langle (m,r),\alpha\rangle$ whose first component belongs to  $\LLn$, and the second component belongs to $L_{p+1}$, i.e.,$(m,r)\in \LLn$ and $\alpha\in L_{p+1}$.

 \subsection*{The Lattice reduct  $\big\langle A^{p+1}_{n+1}; \land,\lor,\langle(n,0),\bfm{0}\rangle,\langle(n,0),\bfm{p}\big\rangle$}
 $\big\langle A^{p+1}_{n+1}; \land,\lor,\langle(n,0),\bfm{0}\rangle,\langle(n,0),\bfm{p}\big\rangle$ is
 the bounded distributive lattice given by the diagram depicted in Figure 1.\smallskip

 ${}$\vspace{1cm}

  \bg{figure}[h]
\bg{center}
  \bg{picture}(140,275)
 \multiput(30,-33)(0,60){2}{\circle*{5}}
 \multiput(30,107)(0,60){2}{\circle*{5}}
 \multiput(210,52)(0,60){2}{\circle*{5}}
 \multiput(210,192)(0,60){2}{\circle*{5}}
 \multiput(-50,-70)(0,60){2}{\circle*{5}}
 \multiput(-50,70)(0,60){2}{\circle*{5}}
 \multiput(-50,-140)(260,460){2}{\circle*{5}}
 \multiput(50,-23)(0,60){2}{\circle*{3}}
 \multiput(50,117)(0,60){2}{\circle*{3}}
 \multiput(210,270)(0,30){2}{\circle*{3}}
 \multiput(-25,-58)(0,60){2}{\circle*{3}}
 \multiput(-25,82)(0,60){2}{\circle*{3}}
 \multiput(12,-41)(0,60){2}{\circle*{3}}
 \multiput(12,98)(0,60){2}{\circle*{3}}
 \multiput(-45,-8)(4,2){8}{\circle*{1}}
 \multiput(-45,-68)(4,2){8}{\circle*{1}}
 \multiput(-45,72)(4,2){8}{\circle*{1}}
 \multiput(-45,132)(4,2){8}{\circle*{1}}
 \multiput(0,92)(4,2){17}{\circle*{1}}
 \multiput(0,152)(4,2){17}{\circle*{1}}
 \multiput(90,55)(4,2){17}{\circle*{1}}
 \multiput(90,195)(4,2){17}{\circle*{1}}
 \multiput(90,-5)(4,2){17}{\circle*{1}}
 \multiput(12,-35)(0,10){20}{\line(0,1){3}}
 \multiput(-25,-53)(0,10){20}{\line(0,1){3}}
 \put(78,30){\makebox(0,0){$\cdots $}}
 \multiput(78,100)(0,5){3}{$\cdot$} 
 \multiput(0,60)(0,5){3}{$\cdot$}
 \put(174,140){\makebox(0,0){$\cdots $}}
 \multiput(-50,-72)(0,10){20}{\line(0,1){3}}
 \multiput(30,-34)(0,10){20}{\line(0,1){3}}
 \multiput(50,-15)(0,10){20}{\line(0,1){3}}
 \multiput(210,55)(0,10){20}{\line(0,1){3}}
 \multiput(101,8)(0,10){20}{\line(0,1){3}}
 \multiput(123,20)(0,10){20}{\line(0,1){3}}
 \multiput(188,45)(0,10){20}{\line(0,1){3}}
 \multiput(144,28)(0,10){20}{\line(0,1){3}}
 \multiput(123,12)(0,60){2}{\circle*{5}}
  \multiput(123,152)(0,60){2}{\circle*{5}}
 \multiput(0,13)(4,2){17}{\circle*{1}}
 \multiput(0,-47)(4,2){17}{\circle*{1}}
 \multiput(90,135)(4,2){17}{\circle*{1}}
 \multiput(101,0)(0,60){2}{\circle*{3}}
 \multiput(101,140)(0,60){2}{\circle*{3}}
 \multiput(-50,-120)(0,30){2}{\circle*{3}}
  \multiput(144,22)(0,60){2}{\circle*{3}}
  \multiput(144,162)(0,60){2}{\circle*{3}}
  \multiput(180,177)(4,2){8}{\circle*{1}}
  \multiput(180,237)(4,2){8}{\circle*{1}}
  \multiput(180,97)(4,2){8}{\circle*{1}}
  \multiput(180,37)(4,2){8}{\circle*{1}}
 \multiput(210,253)(0,5){6}{\circle*{1}}
  \multiput(188,41)(0,60){2}{\circle*{3}}
  \multiput(188,181)(0,60){2}{\circle*{3}}
 \multiput(-50,-100)(0,5){6}{\circle*{1}}
 \multiput(-50,-135)(0,5){6}{\circle*{1}}
 \multiput(210,290)(0,5){6}{\circle*{1}}
 \put(-10,-138){\makebox(0,0){\small$\langle (n,0),\bfm{0}\rangle =\bot$}}
 \put(-80,-120){\makebox(0,0){\small$\langle (n,-k),\bfm{0}\rangle $}}
 \put(-85,-90){\makebox(0,0){\small$\langle (n-1,k),\bfm{0}\rangle $}}
 \put(-85,-70){\makebox(0,0){\small$\langle (n-1,0),\bfm{0}\rangle $}}
 \put(-88,70){\makebox(0,0){\small$\langle (0,0),\bfm{p-1}\rangle $}}
 \put(5,-68){\small\makebox(0,0){\small$\langle (n-1,-k),\bfm{0}\rangle $}}
 \put(40,-50){\small\makebox(0,0){\small$\langle (n-2,k),\bfm{0}\rangle $}}
 \put(62,-40){\small\makebox(0,0){\small$\langle (n-2,0),\bfm{0}\rangle $}} \put(87,-28){\small\makebox(0,0){\small$\langle (n-2,-k),\bfm{0}\rangle $}}
 \put(235,50){\small\makebox(0,0){\small$\langle (0,0),\bfm{0}\rangle $}}
 \put(212,33){\small\makebox(0,0){\small$\langle (0,k),\bfm{0}\rangle $}}
 \put(-63,85){\small\makebox(0,0){\small$\langle (0,k),\bfm{p-1}\rangle $}} \put(-28,98){\small\makebox(0,0){\small$\langle (1,-k),\bfm{p-1}\rangle $}} \put(-8,110){\small\makebox(0,0){\small$\langle (1,0),\bfm{p-1}\rangle $}}
 \put(15,122){\small\makebox(0,0){\small$\langle (1,k),\bfm{p-1}\rangle $}}
 \put(190,15){\small\makebox(0,0){\small$\langle (n-m-1,-k),\bfm{0}\rangle $}}
  \put(160,3){\small\makebox(0,0){\small$\langle (n-m-1,0),\bfm{0}\rangle $}}
 \put(145,-9){\small\makebox(0,0){\small$\langle (n-m-1,k),\bfm{0}\rangle $}}
 \put(56,142){\makebox(0,0){\small$\langle (m,-k),\bfm{p-1}\rangle $}}
 \put(82,154){\makebox(0,0){\small$\langle (m,0),\bfm{p-1}\rangle $}}
 \put(108,165){\makebox(0,0){\small$\langle (m,k),\bfm{p-1}\rangle $}}
  \put(138,185){\makebox(0,0){\small$\langle (n-1,-k),\bfm{p-1}\rangle $}}
 \put(255,194){\makebox(0,0){\small$\langle (n-1,0),\bfm{p-1}\rangle $}}
 \put(245,110){\makebox(0,0){\small$\langle (n-1,0),\bfm{1}\rangle $}}
 \put(245,270){\makebox(0,0){\small$\langle (n-1,k),\bfm{p}\rangle $}}
 \put(240,300){\makebox(0,0){\small$\langle (n,-k),\bfm{p}\rangle$}}
 \put(173,323){\makebox(0,0){\small$\top=\langle (n,0),\bfm{p}\rangle$}}
 \put(-75,-10){\makebox(0,0){\small$\langle (0,0),\bfm{1}\rangle$}}
 \put(70,65){\makebox(0,0){\small$\langle (m,-k),\bfm{1}\rangle $}}
  \put(98,78){\makebox(0,0){\small$\langle (m,0),\bfm{1}\rangle $}}
  \put(120,90){\makebox(0,0){\small$\langle (m,k),\bfm{1}\rangle $}}
  \put(148,109){\makebox(0,0){\small$\langle (n-1,-k),\bfm{1}\rangle $}}
 \put(-15,23){\small\makebox(0,0){\small$\langle (1,-k),\bfm{1}\rangle $}}
 \put(8,34){\small\makebox(0,0){\small$\langle (1,0),\bfm{1}\rangle $}}
  \put(30,45){\small\makebox(0,0){\small$\langle (1,k),\bfm{1}\rangle $}}
  \put(-50,5){\small\makebox(0,0){\small$\langle (0,k),\bfm{1}\rangle $}}
  \put(-74,133){\makebox(0,0){\small$\langle (0,0),\bfm{p}\rangle$}}
 \put(72,205){\makebox(0,0){\small$\langle (m,-k),\bfm{p}\rangle $}}
  \put(98,218){\makebox(0,0){\small$\langle (m,0),\bfm{p}\rangle $}}
  \put(122,230){\makebox(0,0){\small$\langle (m,k),\bfm{p}\rangle $}}
  \put(156,249){\makebox(0,0){\small$\langle (n-1,-k),\bfm{p}\rangle $}}
 \put(245,250){\makebox(0,0){\small$\langle (n-1,0),\bfm{p}\rangle $}}
  \put(-12,163){\small\makebox(0,0){\small$\langle (1,-k),\bfm{p}\rangle $}}
 \put(10,174){\small\makebox(0,0){\small$\langle (1,0),\bfm{p}\rangle $}}
  \put(30,183){\small\makebox(0,0){\small$\langle (1,k),\bfm{p}\rangle $}}
  \put(-46,145){\small\makebox(0,0){\small$\langle (0,k),\bfm{p}\rangle $}}
 \e{picture}
 \end{center}
{}\vspace{4.5cm}

\caption{The lattice reduct of $\aAp$}
\end{figure}

 \no Its lattice order $\le$  can be described as follows:\smallskip

 \no If $\langle(m,r),\alpha\rangle,\langle (k,s),\beta\rangle\in A^{p+1}_{n+1}$, then $\langle(m,r),\alpha\rangle\le\langle (k,s),\beta\rangle$ if and only if either
 \be[o$1)$]
 \item  $\bfm{0}\not=\alpha\le_{\bfm{L_{p+1}}}\beta$,  and $(m,r)\le_{\bfm{L_{n+1}^{\omega}}} (k,s)$, or
 \item $\alpha=\beta=\bfm{0}$, and $(k,s)\ge_{\bfm{L_{n+1}^{\omega}}} (m,r)$, or
 \item $\alpha=\bfm{0}$, $\beta\not=\bfm{0}$ and $(n-1,0)\preccurlyeq(m+k,r+s)$, i.e.,\\
   either $n-1<m+k$, or $n-1=m+k$ and $0\le r+s$.
 \ee

 Observe that for all $(m,r),(k,s)\in\LLn$ and for any $\alpha,\beta\in L_n$, it satisfies:
  \bi
  \item If $\langle(m,r),\alpha\rangle\le\langle(k,s),\beta\rangle$, then $\alpha\le_{\bfm{L_{\bfm{p+1}}}}\beta$.
  \item $\langle(m,r),\bfm{0}\rangle\le\langle(k,s),\bfm{p}\rangle$ if and only if $\langle(k,s),\bfm{0}\rangle\le\langle(m,r),\bfm{p}\rangle$.
  \item $\langle(m,r),\bfm{0}\rangle\le\langle(k,s),\bfm{0}\rangle$ if and only if $\langle(k,s),\bfm{p}\rangle\le\langle(m,r),\bfm{p}\rangle$.
 \item If $(m,r)\ge_{\bfm{L_{n+1}^{\omega}}}(n-1,0)$, then
   \bi
   \item for any $\beta>_{\bfm{L_{p+1}}}\bfm{0}$ $\langle(m,r),\bfm{0}\rangle \le \langle (k,s),\beta\rangle$.

   \item for any  $\beta<_{\bfm{L_{p+1}}}\bfm{p}$ $\langle(k,s),\beta\rangle \le \langle (m,r),\bfm{p}\rangle$
    \ei
 \item $\top=\langle(n,0),\bfm{p}\rangle$  and $\bot=\langle(n,0),\bfm{0}\rangle$ are, respectively,   the greatest and smallest element of $\langle A^{p+1}_{n+1},\le\rangle$.
 \item The induced order in  $\LLn\times\{\bfm{p}\}$, in $L^{\omega}_n\times\{\alpha\}$, for $\bfm{0}<_{\bfm{L_{p+1}}}\alpha<_{\bfm{L_{p+1}}}\bfm{p}$, and in $\LLn\times\{\bfm{0}\}$ is total.
    \ei

 Observe also that the lattice operations are:
   \bi
  \item If $\bfm{0}\notin\{\alpha,\beta\}$, then\\
  $ \langle (m,r),\alpha\rangle \lor \langle(k,s),\beta\rangle =\langle \max\{(m,r),(k,s)\},\max\{\alpha,\beta\}\rangle$.\\
   $\langle (m,r),\alpha\rangle\land \langle(k,s),\beta\rangle =\langle\min\{(m,r),(k,s)\},\min\{\alpha,\beta\}\rangle$.
\item
$\langle (m,r),\bfm{0}\rangle \lor \langle(k,s),\bfm{0}\rangle =
  \langle\min\{(m,r),(k,s)\},\bfm{0}\rangle$.\\
$\langle (m,r),\bfm{0}\rangle\land \langle(k,s),\bfm{0}\rangle =
   \langle\max\{(m,r),(k,s)\},\bfm{0}\rangle$.
\item If $(k,s)\ge_{\bfm{\LLn}}(n-1,0)$, then \\
  $\langle (m,r),\alpha\rangle \lor \langle (k,s),\bfm{0}\rangle =\langle (m,r),\alpha\rangle$ for any $\alpha\not=\bfm{0}$\\
  $\langle (m,r),\alpha\rangle \land \langle (k,s),\bfm{p}\rangle =\langle (m,r),\alpha\rangle$, for any $\alpha\not=\bfm{p}$
 \item If $\alpha\not=\bfm{0}$ and $(k,s)\le_{\bfm{\LLn}}(n-1,0)$, then\\
 $\langle (m,r),\alpha\rangle \lor \langle (k,s),\bfm{0}\rangle =
 \langle (m,r),\alpha\rangle\lor
 \langle(n-1-k,-s), \alpha\rangle.$\\
 $\langle (m,r),\alpha\rangle \land \langle (k,s),\bfm{0}\rangle =
 \langle (n-1-m,-r),\bfm{0}\rangle\land
 \langle(k,s), \bfm{0}\rangle.$
  \ei

  \subsection*{The   partial-ordered  monoid $\langle A^{p}_{n+1};\odot,\top;\le\rangle$}

 For any  $a=\langle(m,r),\alpha\rangle$ and $b=\langle(k,s),\beta\rangle\in A^{p+1}_{n+1}$ we define
 $a\odot b=b\odot a$  as follows
\bi
 \item
 If $\bfm{0}<_{\bfm{L_{p+1}}}\alpha,\beta$ and $\alpha\ast\beta\not=\bfm{0}$, then\\
 $a\odot b=\langle(m,r),\alpha\rangle\odot\langle(k,s),\beta\rangle =\langle(m,r)\ast (k,s),\alpha\ast\beta\rangle$.
 \item  If $\bfm{0}<_{\bfm{L_{p+1}}}\alpha,\beta$ and $\alpha\ast\beta=\bfm{0}$, then\\
  $a\odot b=\langle(m,r),\alpha\rangle\odot\langle(k,s),\beta\rangle =\langle\min\{(n,0),(2n-(m+k+1),-(r+s))\},\bfm{0}\rangle$.\\
  Observe that in this case $\bfm{p}\notin\{\alpha,\beta\}$
  \item If $\alpha\not=\bfm{0}$ and $\beta=\bfm{0}$ then\\
      $a\odot b=\langle(m,r),\alpha\rangle\odot\langle(k,s),\bfm{0}\rangle=
      \langle(m,r)\to (k,s),\bfm{0}\rangle$.
 \item If $\alpha=\bfm{0}$ and $\beta=\bfm{0}$ then\\
  $a\odot b=\langle(m,r),\bfm{0}\rangle\odot\langle(k,s),\bfm{0}\rangle=  \langle\min\{(n,0),(m+k+1,r+s)\},\bfm{0}\rangle$.
\ei
\bg{lemma}\label{L: neu-ab} For any $a\in A^{p+1}_{n+1}$, $a\odot\bot=\bot$, and $a\odot \top=a$.\e{lemma}
{\em Proof:} The statement follows from
 the fact that for any $(m,r)\in \LLn$, we have $(m,r)=(n,0)\ast(m,r)=(n,0)\to (m,r)$. \hfill$\Box$\smallskip

Taking into account that $\bfm{\LLn}$ is an enriched Wajsberg hoop and $\mathds{Z}\stackrel{\to}{\otimes}\mathds{Z}$ is a totally ordered group, it is straightforward, but tedious, to check that $\odot$ is associative (see Appendix A), and so
 $ \langle A^{p+1}_{n+1};\odot,\top\rangle$ is a  commutative
 monoid,  moreover:
 \begin{lemma}\label{L:*monoton3}
  For any $a,b,c\in  A^{p+1}_{n+1}$, $b\le c$ implies $a\odot b\le a\odot c$.
\end{lemma}
{\em Proof:} See Appendix B.\hfill$\Box$\smallskip

Therefore $ \langle A^{p}_{n+1};\odot,\top;\le\rangle$ \emph{is an integral partial-ordered (lattice-ordered) commutative monoid}.

 \subsection*{The residual $\tdiv$}
In what follows   $\neg (m,r)=(m,r)\to (0,0)=(n-m,-r)$ represents the usual negation in $\bfm{\LLn}$, $n>0$. Thus for any $\alpha\in L_{p+1}$, $\neg\, \alpha=\bfm{p}-\alpha$.\smallskip

 For every $\langle(m,r),\alpha\rangle\in A^{p+1}_{n+1}$ we define
 \begin{align*}
 \lefteqn{\sim\!\langle(m,r),\alpha\rangle}\hspace{0.5cm}\\
 &=\left\{\begin{array}{ll}
  \langle(m,r),\neg\alpha\rangle =\langle(m,r),\bfm{p}-\alpha\rangle &\hbox{if } \alpha\in\{\bfm{0},\bfm{p}\} \\
   \langle(n-1-m,-r),\neg\alpha\rangle=\langle(n-1-m,-r),\bfm{p}-\alpha\rangle& \hbox{if } \alpha\notin\{\bfm{0},\bfm{p}\}.
  \end{array} \right.
  \end{align*}
 We are going to prove that we can apply Lemma \ref{L:Monid-invo} to the integral partially ordered commutative monoid $\langle A^{p+1}_{n+1};\odot,\top;\le\rangle$ and $\sim$. First we  prove that $\sim$ is an involution.
  \bg{lemma}
 For any $a,b\in A^{p+1}_{n+1}$ it satisfies
  $a\le b \Lra \,\sim\! b\le \sim\! a$,.  That is,  $\sim$ is an involution.
  \e{lemma}
  {\em Proof:}
  Let $a=\langle(m,r),\alpha\rangle$ and $b=\langle(k,s),\beta\rangle$ be arbitrary elements of $A^{p+1}_{n+1}$ such  that $a\le b$.
  We procedure by cases.
  \be[\bf 1)]
  \item  If $\alpha,\beta\in\{\bfm{0},\bfm{p}\}$,  for $ \alpha=\beta$ the claim is trivial.\\
       If $\alpha=\bfm{0}$ and $\beta=\bfm{p}$, then
  \beqn a\le b
  &\mbox{ iff } & (n-1,0)\preccurlyeq(m+k,r+s)\\
  &\mbox{ iff } &  \sim\!b=\langle(k,s),\bfm{0}\rangle\le \langle(m,r),\bfm{p} \rangle=\sim\!a.
  \eeqn
  \item If $\alpha=\bfm{0}$ and $\beta\notin\{\bfm{0},\bfm{p}\}$, then $\neg \beta\notin\{\bfm{0},\bfm{p}\}$, $(0,0)\preccurlyeq (k,s)\preccurlyeq (n-1)$,   and
  \beqn
  a\le b
  &\mbox{ iff } & (n-1,0)\preccurlyeq(m+k,r+s)\\
  &\mbox{ iff } &  \sim\!b=\langle(n-1-k,-s),\neg \beta\rangle\le \langle(m,r),\bfm{p} \rangle=\sim\!a.
  \eeqn
  \item If $\alpha,\beta\notin\{\bfm{0},\bfm{p}\}$, then $\neg\alpha,\neg\beta\notin\{\bfm{0},\bfm{p}\}$, $(0,0)\preccurlyeq (m,r),(k,s)\preccurlyeq (n-1)$,  and
  \beqn
  a\le b
  &\mbox{ iff } & \alpha\le_{\bfm{L_{p+1}}}\beta \mbox{ and } (m,r)\preccurlyeq(k,s)\\
  &\mbox{ iff } & \neg\beta \le_{\bfm{L_{p+1}}}\neg\alpha  \mbox{ and } (n-1-k,-s)\preccurlyeq(n-1-m,-r)\\
  &\mbox{ iff } &  \sim\!b=\langle(n-1-k,-s),\neg \beta\rangle\le \langle(n-1-m,-r),\neg\alpha \rangle=\sim\!a.
  \eeqn
  \item If $ \alpha\notin\{\bfm{0},\bfm{p}\}$ and $\beta=\bfm{p}$, then $ \neg\alpha\notin\{\bfm{0},\bfm{p}\}$ and
  \beqn
  a\le b &\mbox{ iff }&\bfm{0}<_{\bfm{L_{p+1}}}\alpha<_{\bfm{L_{p+1}}}\bfm{p} \mbox{ and } (m,r)\preccurlyeq(k,s)\\
  &\mbox{ iff } & \bfm{0}=\neg\,\bfm{p}<_{\bfm{L_{p+1}}}\neg \alpha \mbox{ and } (0,0)\preccurlyeq (k-m,s-r)\\
    &\mbox{ iff } & \bfm{0}<_{\bfm{L_{p+1}}}\neg \alpha  \mbox{ and } (n-1,0)\preccurlyeq (k+[(n-1)-m],s-r)\\
  &\mbox{ iff } &  \sim\!b=\langle(k,s), \bfm{0}\rangle\le \langle(n-1-m,-r),\neg\alpha \rangle=\sim\!a.
  \eeqn
 \ee
 So, in all cases the statement is true. \hfill $\Box$ \medskip

 We need to chek that $\odot$ and $\sim$ satisfying the hypothesis of Lemma \ref{L:Monid-invo}.

 \bg{lemma}\label{L:sim} For all $(m,r),(k,s)\in \LLn$, and for any $\beta,\alpha\in L_{p+1}$
 we have
 \[\langle(m,r),\alpha\rangle\odot\langle(k,s),\beta\rangle=\bot\mbox{ iff }\langle(m,r),\alpha\rangle\le \sim\!\langle(k,s), \beta\rangle\]
 \e{lemma}
 {\em Proof:} Recall that $\bot=\langle(n,0),\bfm{0}\rangle$.
 Consider $(m,r),(k,s)$ arbitrary elements of $\LLn$. Observe that for any $\alpha,\beta\in L_{p+1}$,  $\langle(m,r),\alpha\rangle\odot\langle(k,s),\beta\rangle=\bot$ implies $\alpha\ast\beta=\bfm{0}$, which is equivalent to $ \alpha\le_{\bfm{L_{p+1}}}\neg\beta$. We procedure by cases on $\beta$
 \be[\bf 1)]
 \item If \fbox{$\beta=\bfm{p}$}, then
 \beqn
 \langle(m,r),\alpha\rangle\odot\langle(k,s),\bfm{p}\rangle=\bot& \mbox{iff} & \alpha=\bfm{0}=\neg\bfm{p}  \mbox{ and } (k,s)\to(m,r)=(n,0)\\
 & \mbox{iff} & \alpha=\bfm{0}\mbox{ and } (k,s)\le_{\bfm{\LLn}}(m,r).\\
 & \mbox{iff} & \langle(m,r),\alpha\rangle\le\langle(k,s),\bfm{0}\rangle=\sim\langle(k,s),\bfm{p}\rangle \eeqn
 \item If \fbox{$\bfm{0}<_{\bfm{L_{p+1}}}\beta<_{\bfm{L_{p+1}}}\bfm{p}$}, then $\bfm{0}<_{\bfm{L_{p+1}}}\neg\beta<_{\bfm{L_{p+1}}}\bfm{p}$, We have two possibilities,\\ either
 \underline{$\alpha=\bfm{0}$}, and so, since $\sim$ is an involution,
 \beqn
 \langle(m,r),\bfm{0}\rangle\odot\langle(k,s),\beta\rangle=\bot  & \mbox{iff }&  (k,s)\to(m,r)=(n,0) \\
 & \mbox{iff } & (k,s)\le_{\LLn}(m,r)\\
 & \mbox{iff } & \langle(k,s),\beta\rangle\le\langle(m,r),\bfm{p}\rangle =\sim \langle(m,r),\bfm{0}\rangle.
 \eeqn
 or
 \underline{$\alpha\not=\bfm{0}$}, in this case $\bfm{0}<_{\bfm{L_{p+1}}}\alpha<_{\bfm{L_{p+1}}}\bfm{p}$, then
 \beqn
 \lefteqn{\langle(m,r),\alpha\rangle\odot\langle(k,s),\beta\rangle=\bot}\hspace{1.75cm}\\
 & \mbox{iff }& \alpha\le_{\bfm{L_{p+1}}}\neg\beta \mbox{ and } (n,0)\preccurlyeq \big(2n-(m+k+1),-(r+s)\big) \\
 & \mbox{iff } &\alpha\le_{\bfm{L_{p+1}}}\neg\beta \mbox{ and } (m,r)\le_{\LLn}(n-1-k,-s)\\
 & \mbox{iff } & \langle(m,r),\alpha\rangle\le\langle(n-k-1,-s),\neg\beta\rangle =\sim \langle(k,s),\beta\rangle.
 \eeqn

 \item If \fbox{$\beta=\bfm{0}$}.\\
 Either  \underline{$\alpha=\bfm{p}$}, and so, by case \textbf{1)} and commutativity of $\odot$
 \beqn
  \langle(m,r),\bfm{p}\rangle\odot\langle(k,s),\bfm{0}\rangle=\bot& \mbox{iff} &
  \langle(k,s),\bfm{0}\rangle\le\sim\langle(m,r),\bfm{0}\rangle\\
  & \mbox{iff} &  \langle(m,r),\bfm{p}\rangle\le\langle(k,s),\bfm{p}\rangle =\sim \langle(k,s),\bfm{0}\rangle,
 \eeqn
 or  \underline{$\bfm{0}<_{\bfm{L_{p+1}}}\alpha<_{\bfm{L_{p+1}}}\bfm{p}$}, and the property follows from  then  involutive property of $\sim$ and the case 2).\\
 or  \underline{$\alpha=\bfm{0}$}, and so
 \beqn
  \langle(m,r),\bfm{0}\rangle\odot\langle(k,s),\bfm{0}\rangle=\bot & \mbox{iff} & (n,0)\preccurlyeq(m+k+1,r+s)\\
  &\mbox{iff} & (n-1,0)\preccurlyeq(m+k,r+s)\\
 & \mbox{iff} &\langle(m,r),\bfm{0}\rangle< \langle(k,s),\bfm{p}\rangle=\,\sim\! \langle(k,s),\bfm{0}\rangle\ \ \ \hfill \Box
 \eeqn
 \ee
 If we take
 $x\tdiv y=\sim(x\odot \sim y),$
 then, from Lemmata \ref{L:sim} and \ref{L:Monid-invo}, we deduce
 \bg{theorem}  For any $a, b, c\in A^{p+1}_{n+1}$,
 $a\odot b\le c \mbox{ if and only if } b\le a\tdiv c.$
 And so $\bfm{A^{p+1}_{n+1}}$ is a bounded integral commutative residuated lattice.\hfill$\Box$
 \e{theorem}

The operation  $\tdiv$ can be described as follows.
Let $(m,r),(k,s)$ be arbitrary elements of $\LLn$, then
 \beqn
 \langle(m,r),\bfm{p}\rangle\tdiv \langle(k,s),\bfm{p}\rangle &=& \langle(k,s),\bfm{0}\rangle\tdiv \langle(m,r),\bfm{0}\rangle= \langle(m,r)\to(k,s),\bfm{p}\rangle,\\
 \langle(m,r),\bfm{p}\rangle\tdiv \langle(k,s),\bfm{0}\rangle
 &=& \langle(k,s),\bfm{p}\rangle\tdiv \langle(m,r),\bfm{0}\rangle=\langle(m,r)\ast(k,s),\bfm{0}\rangle,\\
\langle(m,r),\bfm{0}\rangle\tdiv \langle(k,s),\bfm{p}\rangle
 &=& \langle(k,s),\bfm{0}\rangle\tdiv \langle(m,r),\bfm{p}\rangle\\
 &=&\langle\min\{(n,0),(m+k+1,r+s)\},\bfm{p}\rangle.
 \eeqn
If $\bfm{0}<_{\bfm{L_{p+1}}}\alpha\le_{\bfm{L_{p+1}}} \beta\le_{\bfm{L_{p+1}}} \bfm{p}$, then $\alpha\ast\neg\beta =\bfm{0}$, and
\beqn
 \langle(m,r),\alpha\rangle\tdiv \langle(k,s),\beta\rangle
 &=& \langle(m,r)\to(k,s),\bfm{p}\rangle
 =\langle(m,r)\to(k,s),\alpha\to \beta\rangle.
 \eeqn
 If $\bfm{0}<_{\bfm{L_{p+1}}}\alpha\not\le_{\bfm{L_{p+1}}}\beta$ and $ \beta\not=\bfm{0} $, then $\beta\not=\bfm{p}$ and, since $\alpha\to\beta=\neg(\alpha\ast\neg\beta)\not=\bfm{p}$, we have
\beqn
 \langle(m,r),\alpha\rangle\tdiv \langle(k,s),\beta\rangle
 &=& \langle(m,r)\to(k,s),\alpha\to\beta\rangle.
 \eeqn
 If $\bfm{0}<_{\bfm{L_{p+1}}}\alpha<_{\bfm{L_{p+1}}}\bfm{p}$, then $\bfm{0}<_{\bfm{L_{p+1}}}\neg\alpha<_{\bfm{L_{p+1}}}\bfm{p}$, and  taking into account that $(m,r)\ast (k,s)\le_{\bfm{\LLn}}(m,r)\le_{\bfm{\LLn}}(n-1,0)$, we obtain
\beqn
 \langle(m,r),\alpha\rangle\tdiv \langle(k,s),\bfm{0}\rangle
 &=  & \langle\min{(n-1,0),(m+k+1,r+s)},\neg\,\alpha\rangle.
 \eeqn

  \bg{remark} Observe that for any $a\in \Ap$, $\sim a=a\tdiv \langle (n,0),\bfm{0}\rangle= a\tdiv \bot$.\e{remark}

 \section{The Involutive bounded commutative re\-si\-dua\-ted lattice $\aAp$}

  Summarizing  the above results, we have

 \bg{theorem} For any integers $n,p>0$,  the algebra
  \[\aAp=\left\langle \Ap ;\odot,\tdiv,\land,\lor,\langle (n,0),\bfm{0}\rangle,\langle (n,0),\bfm{p}\rangle\right\rangle\]
  is a \brl, i.e.,it satisfies:
  \bi
   \item $\left\langle \Ap; \odot ,\langle (n,0),\bfm{p}\rangle\right\rangle$ is a commutative monoid,
 \item $\left\langle \Ap; \lor, \land,\langle (n,0),\bfm{0}\rangle,\langle (n,0),\bfm{p}\rangle\right\rangle$ is a bounded distributive\footnote{Because it  contains no copy of pentagon and diamond.}
 lattice with\\ smallest element $\bot=\langle (n,0),\bfm{0}\rangle$, and greatest element $\top=\langle (n,0),\bfm{p}\rangle$
 \item  It satisfies  residuation
 condition:  for any $a,b,c\in \Ap$
 \bg{align}\label{Eq:resd-cond}
a \odot  b \le  c &
 \mbox{ if and only if } a \le b\tdiv c,\end{align}
 where  $\le $ is the order given by the
 lattice structure.
 \ei
  In addition, $\aAp$ is involutive:
  \bi
  \item For any $a\in \Ap$, $\sim\sim a=\big[ x\tdiv \langle (n,0),\bfm{0}\rangle\big]\tdiv \langle (n,0),\bfm{0}\rangle= a$.\hfill $\Box$
  \ei
  \e{theorem}


  \subsection{Subalgebras}

  It is clear that  for any $n,p>0$  the two-element Boolean algebra  $\bfm{L_2}$ and the Chang's algebra $\bfm{L_2^{\omega}}$ are  copies of subalgebras of $\aAp$.
  Moreover,
 \[\widehat{L_{n+1}^{p+ 1}} = \{\left\langle (m,0),\bfm{q}\right\rangle:  \left\langle (m,0),\bfm{q}\right\rangle\in \Ap\}
    =\{\langle \bfm{m},\bfm{q}\rangle:  \langle \bfm{m},\bfm{q}\rangle\in \Ap \}
 \]
 is universe of a subalgebra of $\aAp$, namely $\widehat{\bfm{L}_{n+1}^{p+ 1}}$ (see Figure 2).
 {}\vspace{0.5cm}

   \bg{figure}[h]
\bg{center}
  \bg{picture}(100,150)
 \multiput(-10,-56)(0,30){2}{\circle*{3}}
 \multiput(-10,24)(0,30){2}{\circle*{3}}
 \multiput(110,14)(0,30){2}{\circle*{3}}
 \multiput(110,93)(0,30){2}{\circle*{3}}
 \multiput(-50,-80)(0,30){2}{\circle*{3}}
 \multiput(-50,0)(0,30){2}{\circle*{3}}
 \multiput(-50,-116)(160,275){2}{\circle*{3}}
 \multiput(50,60)(0,30){2}{\circle*{3}}
  \multiput(50,-20)(0,30){2}{\circle*{3}}
 \multiput(-50,-80)(0,30){2}{\line(5,3){50}}
 \multiput(33,-30)(0,30){2}{\line(5,3){30}}
 \multiput(33,50)(0,30){2}{\line(5,3){30}}
 \multiput(-50,0)(0,30){2}{\line(5,3){50}}
 \multiput(90,2)(0,30){2}{\line(5,3){20}}
 \multiput(90,80)(0,30){2}{\line(5,3){20}}
  \multiput(-50,-115)(160,195){2}{\line(0,1){75}}
  \multiput(-10,-55)(0,70){2}{\line(0,1){40}}
  \multiput(50,-20)(0,70){2}{\line(0,1){40}}
  \multiput(-50,-10)(160,22){2}{\line(0,1){40}}
 \multiput(10,64)(5,3){3}{$\cdot$}
 \multiput(10,34)(5,3){3}{$\cdot$}
 \multiput(70,98)(5,3){3}{$\cdot$} 
 \multiput(70,68)(5,3){3}{$\cdot$} 
 \multiput(70,18)(5,3){3}{$\cdot$} 
 \multiput(70,-12)(5,3){3}{$\cdot$} 
 \multiput(-50,-34)(0,7){3}{\line(0,1){3}}
 \multiput(-10,-10)(0,7){3}{\line(0,1){3}}
 \multiput(50,25)(0,7){3}{\line(0,1){3}}
 \multiput(110,56)(0,7){3}{\line(0,1){3}}
 \multiput(10,-15)(5,3){3}{\circle*{1}}
 \multiput(10,-45)(5,3){3}{\circle*{1}}
 \put(-20,-118){\makebox(0,0){\small$\langle \bfm{n},\bfm{0}\rangle =\bot$}}
 \put(-80,-80){\makebox(0,0){\small$\langle \bfm{n-1},\bfm{0}\rangle $}}
 \put(-78,0){\makebox(0,0){\small$\langle \bfm{0},\bfm{p-1}\rangle $}}
 \put(18,-63){\small\makebox(0,0){\small$\langle \bfm{n-2},\bfm{0}\rangle $}}
 \put(127,12){\small\makebox(0,0){\small$\langle \bfm{0},\bfm{0}\rangle $}}
 \put(15,15){\small\makebox(0,0){\small$\langle \bfm{1},\bfm{p-1}\rangle $}}
  \put(90,-28){\small\makebox(0,0){\small$\langle \bfm{n-m-1},\bfm{0}\rangle $}}
 \put(77,54){\makebox(0,0){\small$\langle \bfm{m},\bfm{p-1}\rangle $}}
 \put(148,94){\makebox(0,0){\small$\langle \bfm{n-1},\bfm{p-1}\rangle $}}
 \put(138,44){\makebox(0,0){\small$\langle \bfm{n-1},\bfm{1}\rangle $}}
 \put(140,160){\makebox(0,0){\small$\top=\langle \bfm{n},\bfm{p}\rangle$}}
 \put(-70,-50){\makebox(0,0){\small$\langle \bfm{0},\bfm{1}\rangle$}}
  \put(68,7){\makebox(0,0){\small$\langle \bfm{m},\bfm{1}\rangle $}}
 \put(6,-34){\small\makebox(0,0){\small$\langle \bfm{1},\bfm{1}\rangle $}}
  \put(-70,30){\makebox(0,0){\small$\langle \bfm{0},\bfm{p}\rangle$}}
  \put(68,87){\makebox(0,0){\small$\langle \bfm{m},\bfm{p}\rangle $}}
 \put(135,118){\makebox(0,0){\small$\langle \bfm{n-1},\bfm{p}\rangle $}}
 \put(5,50){\small\makebox(0,0){\small$\langle \bfm{1},\bfm{p}\rangle $}}
 \put(150,-85){\small\makebox(0,0){\small\fbox{$\langle \bfm{m},\alpha\rangle =\langle (m,0),\alpha\rangle$}}}
 \e{picture}
 \end{center}
{}\vspace{3.5cm}

\caption{The lattice reduct of $\widehat{\bfm{L}_{n+1}^{p+1}}$}
\end{figure}

 Notice that   $\widehat{\bfm{L}_{n+1}^{2}}$ is the algebra $\widehat{\bfm{L_{n+1}}}$, defined in \cite[page 1131]{CT12} to give an example of an algebra bellowing to  $\WLn\smsm\mbb{WL_\mbb{n}}$  having Boolean retraction term. If $q>0$ and $q$ divides $p$ or $q=1$, then $\bfm{A^{q+1}_{n+1}}$ and $\widehat{\bfm{L}_{n+1}^{q+ 1}}$ are  copies of   subalgebras of $\aAp$.
 \smallskip

 It is straightforward that, up isomorphism, there are no more proper subalgebras than those described above. So, if $S$ represents de class operator \emph{isomorphic subalgebras} and $I$ represents the operator isomorphic images, then for any $n,p>0$
  \bg{align*}S(\aAp)= & I(\{\aAp,\bfm{L_2},\bfm{L_2^\omega},\widehat{\bfm{L}_{n+1}^{p+ 1}},\bfm{A^{2}_{n+1}},\widehat{\bfm{L}_{n+1}^{2}}\}\cup\\
  & \{\bfm{A^{q+1}_{n+1}}: q \mbox{ divides }p\} \cup\{\widehat{\bfm{L}_{n+1}^{q+ 1}}: q \mbox{ divides }p\})
  \end{align*}


  \subsection{Implicative filters and homomorphic images}

  The variety $\BRL$ of \brl s is $\top$-regular and  the set of all implicative filters is isomorphic to the set of congruence relations, both ordered by inclusion (see \cite{GJKO} and \cite{CT12} and the references given therein, see also \cite{T23}). So, we will see what the implicative filters of $\aAp$ are.\smallskip

  We recall that in a \brl\ $\A=\langle A; \to, \ast, \land,\lor,\bot,\top\rangle$, a subset $F$ of $A$ is an  \emph{implicative filter} provided $\top\in F$,   $F$ is closed by the monoidal operation $\ast$ and  it is  increasing, or equivalently, $\top\in F$  and $x,x\to y\in F$ implies $y\in F$. Then $\theta(F)=\{(x,y)\in A\times A: (x\to y)\ast (y\to x)\in F\}$ is a congruence relation such that $\top/\theta(F)=F$; and for any congruence relation $\theta$, $\top/\theta$ is an implicative filter such that $\theta(\top/\theta)=\theta$.
  \smallskip

  For any $n,p>0$  consider the following subsets of $\Ap$:
  \be[({i}f1)]
  \item $\{\top\}$
  \item $f_\omega=\{\langle(n,r),\bfm{p}\rangle\in A^{p+1}_{n+1}:r\le 0\}$
  \item $\uparrow\!\langle\bfm{0},\bfm{p}\rangle=\{a\in A^{p+1}_{n+1}:a\ge\langle\bfm{0},\bfm{p}\rangle\}$.
      \ee
  Observe $\{\top\}\subseteq f_\omega\subseteq\uparrow\!\langle\bfm{0},\bfm{p}\rangle$.\smallskip

  Notice  that $\langle\bfm{n-1},\bfm{p-1}\rangle =\max\{a\in \Ap: a\not\ge \langle\bfm{0},\bfm{p}\rangle\}$, then $a\not\ge \langle\bfm{0},\bfm{p}\rangle$   if and only if  $a\le \langle\bfm{n-1},\bfm{p-1}\rangle$, and so $\Ap=\uparrow\!\langle\bfm{0},\bfm{p}\rangle\cup\! \downarrow\!\!\langle\bfm{n-1},\bfm{p-1}\rangle$ (disjoint union), where $\downarrow\!\!\langle\bfm{n-1},\bfm{p-1}\rangle=\{a\in \Ap: a\le\langle\bfm{n-1},\bfm{p-1}\rangle \}$. \smallskip

  Since  for any $k>n,p$ it satisfies
   \bi
   \item if $(m,s)\in L_{n+1}^\omega\smsm\{(n,r): r\le 0 \}$, then $(m,s)^{k}=(0,0)$,
   \item $\langle\bfm{0},\bfm{p}\rangle^k=\langle\bfm{0},\bfm{p}\rangle$, and
   \item $\langle\bfm{n-1},\bfm{p-1}\rangle^k=\langle\bfm{n},\bfm{0}\rangle=\bot$,
  \ei
   then we have
  \bg{lemma}\label{L:Local} For all $n,p>0$  $\aAp$  has only three proper implicative filters:  namely $\uparrow\!\langle\bfm{0},\bfm{p}\rangle$, $f_\omega$ and $\{\top\}$.\hfill$\Box$
  \e{lemma}
  And so
\bg{corollary} For any $n,p>0$
 $\aAp$ is local and its proper maximal implicative filter is
 $
 Rad(\aAp) =\,\uparrow\!\langle\bfm{0},\bfm{p}\rangle$\!
 \footnote{For any \brl\ $\A$ its radical, denoted by $ Rad(\A)$, is the intersection of its maximal implicative filters. In our case the radical is equal to the unique maximal implicative filter.}.\hfill $\Box$
 \e{corollary}

   Therefore, for any $n,p>0$ $\aAp$ has  three non-trivial homomorphic images, up to isomorphism:  $\aAp/\{\top\}\cong \aAp$, $\aAp/f_\omega\cong \widehat{\bfm{L}_{n+1}^{p+1}}$ and $\aAp/\!\uparrow\!\!\langle\bfm{0},\bfm{p}\rangle\cong \bfm{L_p}$. Hence, if $H$ represents the
class operator \emph{homomorphic images}, we have $H(\aAp)=I(\{\aAp,\widehat{\bfm{L}_{n+1}^{p+ 1}},\bfm{L_p}\}\cup\{trivial\ algebra\})$.

\subsection{Terms in $\aAp$ and General Apple Property(GAP)}

 We are going to show that for any $n,p>0$, the term $t_{np}(x)= (n+1).x^{\max\{n+1,p\}}$ is boolean term  for $\aAp$, and so $t_{np}(x)$ is also boolean term for its generated variety $HSP(\aAp)$. Thus,  by  \cite[Theorem 3.10]{T23}, $HSP(\aAp)$ has the General Apple Property (GAP). We will also see that  the class  $HSP(\aAp)_\mbb{\,SS}$ of all semisimple members of $HSP(\aAp)$ is a variety, in fact it is contained in the variety $\EMp$, the subvariety of $\BRL$ given by the equation $x\lor \neg x^p\approx \top$. Therefore there is $k>0$ such that $k.x^k$ is boolean and radical term for $HSP(\aAp)$. First we prove a previous result.

\bg{theorem}\label{T:ApLocal} For any  $n,p>0$  and any   $a\in A^{p+1}_{n+1}$ it satisfies:
   \be[$(a)$]
   \item $a\lor\!\sim\! a^p\in Rad(\aAp)$
   \item If $p\le n$, then
  \begin{equation}a\not\ge \langle\bfm{0},\bfm{p}\rangle \mbox{ if and only if }a^{n+1}=\bot.\label{e:a}\end{equation}
  \item If $n<p$, then $\langle\bfm{n-1},\bfm{p-1}\rangle^{p-1}=\, \sim\!\!\langle\bfm{n-1},\bfm{p-1}\rangle$, i.e $\langle\bfm{n-1},\bfm{p-1}\rangle$  is cyclic element of $\aAp$ (cf. \cite{T94}), and
  \begin{equation}a\not\ge \langle\bfm{0},\bfm{p}\rangle \mbox{ if and only if }a^{p}=\bot.\label{e:b}\end{equation}
  \ee
 \e{theorem}
 {\em Proof of $(a)$:} Assume $a\notin Rad(\aAp)$, then $a\le \langle\bfm{n-1},\bfm{p-1}\rangle$, and so there are $(m,r)\in L^\omega_{n+1}$ and $\alpha<\bfm{p}$ such that $a=\langle(m,r),\alpha\rangle$ then, since $\alpha^p=\bfm{0}$, there is $(k,s)\in L^\omega_{n+1}$ such that $a^p=\langle(k,s),\bfm{0}\rangle$. Thus $\sim\! a^p =\langle(k,s),\bfm{p}\rangle\in Rad(\aAp)$. This implies $a\lor\sim a^p\in Rad(\aAp)$.\smallskip

 \no On the other hand,  observe that the converse implication of (\ref{e:a}) and of (\ref{e:b}) follows from the fact that for any $k>0$   $a^k=\bot$ implies $a\notin \uparrow\!\!\langle\bfm{0},\bfm{p}\rangle$. 

 \no To see the direct implication of (\ref{e:a}) we suppose  $p\le n$. It suffices to see  $\langle\bfm{n-1},\bfm{p-1}\rangle^{n+1}=\bot$.
 \beqn
 \lefteqn{\langle\bfm{n-1},\bfm{p-1}\rangle^{n+1}}\hspace{1cm}\\
 &=&[\langle\bfm{n-1},\bfm{p-1}\rangle^{p-1}\odot \langle\bfm{n-1},\bfm{p-1}\rangle]\odot\langle\bfm{n-1},\bfm{p-1}\rangle^{n-p+1}\\
 &=&[\langle\bfm{n}-(\bfm{p}-\bfm{1}),\bfm{1}\rangle\odot \langle\bfm{n-1},\bfm{p-1}\rangle]\odot\langle\bfm{n-1},\bfm{p-1}\rangle^{n-p+1}\\
 &=&\langle\bfm{p-1},\bfm{0}\rangle \odot\langle\bfm{n-1},\bfm{p-1}\rangle^{n-p+1}\\
 &=&\langle(\bfm{n-1})^{n-p+1}\to(\bfm{p-1}),\bfm{0}\rangle\\&=&
 \langle(\min\{\bfm{n},([n-p+1](n-(n-1))+(p-1),0)\},\bfm{0}\rangle\\&=&
 \langle(\min\{\bfm{n},((n-p+1)+(p-1),0)\},\bfm{0}\rangle\\&=&
 \langle\bfm{n},\bfm{0}\rangle =\bot,
 \eeqn
 this closes the proof of $(a)$.

 \no To show $(b)$ suppose  $p> n$, then there are $m>0$ and $0< q<n$ such that $p=m\cdot n+q$, so
  \beqn
\langle\bfm{n-1},\bfm{p-1}\rangle^{p-1}
 &=&\langle\bfm{n-1},\bfm{p-1}\rangle^{m\cdot n}\odot\langle\bfm{n-1},\bfm{p-1}\rangle^{q-1}\\
 &=& \langle\bfm{n-1}^{m\cdot n},\bfm{p-1}^{m\cdot n}\rangle\odot\langle\bfm{n-1}^{q-1},\bfm{p-1}^{q-1}\rangle\\
 &=& \langle\bfm{0},\bfm{q}\rangle\odot\langle\bfm{n- (q-1)},\bfm{p-(q-1)}\rangle\\
 &=& \langle\bfm{0},\bfm{1}\rangle= \sim\!\!\langle\bfm{n-1},\bfm{p-1}\rangle
 \eeqn
 and so $\langle\bfm{n-1},\bfm{p-1}\rangle$ is cyclic.\\
 Moreover,
 since $a\not\ge \langle\bfm{0},\bfm{p}\rangle$ if and only if  $a\le \langle\bfm{n-1},\bfm{p-1}\rangle$, by the above, we have that for any $a\not\ge \langle\bfm{0},\bfm{p}\rangle$,
  \[a^{p}\le \langle\bfm{n-1},\bfm{p-1}\rangle^{p}= \langle\bfm{n-1},\bfm{p-1}\rangle^{p-1}\odot\langle\bfm{n-1},\bfm{p-1}\rangle =\bot.\]
  That shows the direct implication of (\ref{e:b}), and the proof of $(b)$ is complete.\hfill $\Box$
  \smallskip

  From the proof of the above theorem we deduce:
  \bg{corollary}
  For any  $n,p>0$ $\aAp $  satisfies:
  \be[$(a)$]
  \item If $k\ge\max\{n+1,p\}$ and $a\in \Ap$, then $a^k=\bot$ if and only if $a\not\ge\langle\bfm{0},\bfm{p}\rangle$
  \item If $0<k<\max\{n+1,p\}$, then $\langle\bfm{n-1},\bfm{p-1}\rangle^k\not=\bot$
  \ee
  \e{corollary}
  {\em Proof:} Item $(a)$ is an immediate consequence of the preceding theorem.

  \no To see $(b)$ take $0<k<\max\{n+1,p\}$. If $p\le n$, then
  \[\langle\bfm{n-1},\bfm{p-1}\rangle^k\ge \langle\bfm{n-1},\bfm{p-1}\rangle^n= \langle\bfm{n-1},\bfm{0}\rangle>\bot,\]
  and if $p>n$, then
  $\langle\bfm{n-1},\bfm{p-1}\rangle^k\ge\langle\bfm{n-1},\bfm{p-1}\rangle^{p-1}=
  \langle\bfm{0},\bfm{1}\rangle>\bot$.\hfill$\Box$
  \smallskip

 We are going to see that for any $p,n>0$, $\aAp$ has Boolean and radical term
  \bg{theorem}\label{T:Bterm} For any $p,n>0$, the following are satisfied:
  \be[$(a)$]
  \item $t_{np}(x)=(n+1).x^{\max\{n+1,p\}}$ is Boolean  term for $\aAp$, that is

 $\aAp\models t_{np}(x)\lor \neg t_{np}(x)\approx\top$
\item $
  Rad(\aAp)  =\{a\in \Ap:t_{np}^{\aAp}(a)=\top\} $\\
  ${}$ \hspace{2.1cm}$=\{a\in \Ap:\mbox{ for any } k>0\, (n+1).a^{k}=\top\}$

 \item For any $k>0$, $\aAp\models(n+1).(x\lor\neg x^p)^k\approx\top$.
 \ee
 \e{theorem}{\em Proof:}
 Take $m=\max\{n,p+1\}$.
 From Theorem \ref{T:ApLocal}, we deduce that  any $a\notin Rad(\aAp))= \uparrow\!\!\langle\bfm{0},\bfm{p}\rangle$  satisfies $(n+1).a^m=\bot$,  and so  $\sim\!\!\big((n+1).a^m\big)=\top $.\smallskip

 \no Assume $a=\langle(l,r),\bfm{p}\rangle\in \uparrow\!\langle\bfm{0},\bfm{p}\rangle$, then, since $\uparrow\!\langle\bfm{0},\bfm{p}\rangle$ is an implicative filter, for any $k>0$ $a^k\in\uparrow\!\langle\bfm{0},\bfm{p}\rangle$. Moreover,
 \beqn
 (n+1).a&=&\sim\!\Big(\big(\sim\!\langle(l,r),\bfm{p}\rangle\big)^{n+1}\Big)=\,
 \sim\!\big(\langle(l,r),\bfm{0}\rangle^{n+1}\big)\\&=&
 \sim\!\big(\langle\min\{(n,0),\big((n+1)l+n,nr\big)\},\bfm{0}\rangle\big)\\&=&
 \sim\!\big(\langle\min(n,0),\bfm{0}\rangle\big) = \langle\min(n,0),\bfm{p}\rangle=\top.
 \eeqn
 This proves that $(n+1).x^m$ is both Boolean and radical for $\aAp$, and so item $(a)$ and  the first equality of $(b)$ hold. Moreover, by Theorem 2.12 of \cite{T23}, $\aAp$ is $(n+1)$-radical, and so   the second equality of $(b)$ also holds.\\
  To see $(c)$ observe that by  item $(a)$ of Theorem \ref{T:ApLocal}  for any $a\in \Ap$ $a\lor\! \sim\! a^p\in Rad(\aAp)$. Thus, by item $(b)$, for any $k>0$ the equation $(n+1).(x\lor\neg x^p)^k\approx\top$ holds in $\aAp$. \hfill $\Box$\smallskip

 From  the above theorem, and Teorem 3.17 of \cite{T23} we  deduce
 \begin{corollary}
 For any $n,p>0$ it satisfies
 \be[$(a)$]
   \item $\max\{n+1,p\}.x^{\max\{n+1, p\}}$ is Boolean and radical term for $\aAp$.
   \item $\max\{n+1,p\}.x^{\max\{n+1, p\}}$ is Boolean term for $HSP(\aAp)$.
   \item $HSP(\aAp)_\mbb{\,SS}\subseteq \EMp$ and so $HSP(\aAp)_\mbb{\,SS}$ is a variety.
   \item There is $k>0$ such that $m.x^k$ is boolean and radical term for $HSP(\aAp)$.
 \ee
   \end{corollary}

 Then, by Corollary 2.7 and Theorem  3.10 of \cite{T23},  \emph{$\aAp$ has the General Apple Pro\-perty} (GAP), and so, \emph{its generated variety $HSP(\aAp)$ also has} GAP.
 Therefore, by \cite[Lemma 3.15]{T23}, if $k\ge\max\{n+1,p\}$ then $\aAp\in \mbb{WL}_\mbb{k}$, i.e, $\aAp\models k.x\lor k.(\neg x)=\top$, and so $HSP(\aAp)\subseteq \mbb{WL}_\mbb{k}$. However, from Theorem 5.13 of \cite{CT12} it deduces that the subalgebra $\widehat{\bfm{L}_{n+1}^{2}}$ of $\aAp$  does not belongs to $\WL_\mbb{n}$, and so $\aAp\notin \mbb{WL}_{\mbb{n}}$. Hence we have
  that for any $n,p>0$, if $k=\max\{n+1,p\}$, then  $\aAp\in \mbb{WL}_{\mbb{k}}\smsm \mbb{WL}_{\mbb{n}}$.

\section*{Appendix}
\appendix
\section{Checking the associative of $\odot$ in  $\bfm{A^{p+1}_{n+1}}$}

To  check the associativity of $\odot$ in $ A^{p+1}_{n+1}$,  we consider
 $\langle(m,r),\alpha\rangle,\langle(k,s),\beta\rangle$ and $\langle(l,t),\gamma\rangle$  arbitrary elements of $A^{p+1}_{n+1}$.
 We procedure by cases:

 \subsection{$\alpha=\beta=\gamma=\bfm{p}$\label{A.1}}
 The associativity of $\odot$ follows form the associativity of $\ast$ in $\bfm{\LLn}$.
 \vspace{-0.5cm}

 \subsection{ $\alpha=\bfm{p}$  and $\bfm{0}<_{\bfm{L_{p+1}}}\beta,\gamma$\label{A.2}}
  \subsubsection{ $\bfm{0}=\beta\ast\gamma$\label{A.2.1}}
 In this case  $\bfm{p}\notin\{\beta,\gamma\}$, and
 \beqn
  \lefteqn{\langle(m,r),\bfm{p}\rangle\odot\big[\langle(k,s), \beta\rangle\odot\langle(l,t),\gamma  \rangle\big]}\hspace{1cm}\\
 &=&\langle(m,r),\bfm{p}\rangle\odot\big\langle\min\{(n,0),(2n-(k+l+1),-(s+t))\},
 \bfm{0}  \rangle\\
 &=& \left\{\begin{array}{ll}
      \langle(m,r),\bfm{p}\rangle\odot \langle(n,0),\bfm{0}\rangle\, \ \hbox{if } (n,0)\preccurlyeq (2n-(k+l+1),-(s+t))\\
      \langle(m,r),\bfm{p}\rangle\odot \langle(2n-(k+l+1),-(s+t)),\bfm{0}\rangle, \quad  \hbox{otherwise}
       \end{array}
     \right.\\
 &=& \left\{\begin{array}{ll}
      \langle(m,r)\to(n,0),\bfm{0}\rangle, \quad \ \hbox{if } (k+l+1,s+t)\preccurlyeq (n,0)\\
      \langle(m,r)\to (2n-(k+l+1),-(s+t)),\bfm{0} \rangle, \quad\hbox{otherwise}
       \end{array}
     \right.\\
 &=& \left\{\begin{array}{ll}
      \langle(n,0),\bfm{0}\rangle , \qquad\ \qquad\ \hbox{if } (k+l+1,s+t)\preccurlyeq (n,0)\\
      \langle \min\{(n,0),\big(3n-(m+k+l+1),-(r+s+t)\big)\},\bfm{0} \rangle, \ \hbox{otherwise},
       \end{array}
     \right.\\
 &&\mbox{since }(m,r)\preccurlyeq(n,0)\mbox{ and } (k+l+1,-(s+t))\preccurlyeq (n,0) \mbox{ imply}\\
     && \big(3n-(m+k+l+1),-(r+s+t)\big)\succcurlyeq (n,0), \mbox{ we have}\\
     &=&\langle \min\{(n,0),\big(3n-(m+k+l+1),-(r+s+t)\big)\}, \bfm{0} \rangle.\\
 \eeqn
 ${}$\vspace{-0.95cm}

 \no Moreover,
 \beqn
  \lefteqn{
  \big[\langle(m,r),\bfm{p}\rangle\odot\langle(k,s), \beta\rangle\big]\odot\langle(l,t),\gamma  \rangle}\hspace{1.cm}\\
  &=&\langle(m,r)\ast(k,s),
 \beta  \rangle\odot\langle(l,t),\gamma\rangle\\
  &=&\langle\max\{(0,0),(m+k-n,r+s\},
 \beta  \rangle\odot\langle(l,t),\gamma\rangle\\
&=& \left\{\begin{array}{ll}
      \langle(0,0),\beta\rangle\odot \langle(l,t),\gamma\rangle, \ \hbox{if } (m+k-n,r+s)\preccurlyeq (0,0)\\
      \langle(m+k-n,r+s),\beta\rangle\odot \langle(l+1,t),\gamma\rangle,\quad  \hbox{otherwise}
       \end{array}
     \right.\\&=& \left\{\begin{array}{ll}
      \langle\min\{(n,0),\big(2n-(0+l+1),-t\big),\bfm{0}\rangle,  \hbox{if } (m+k-n,r+s)\preccurlyeq (0,0)\\
      \langle(m+k-n,r+s),\beta\rangle\odot \langle(l+1,t),\gamma\rangle,\quad  \hbox{otherwise}
       \end{array}
     \right.\\
 &=&\left\{\begin{array}{ll}
     \langle(n,0),\bfm{0}\rangle, \qquad \hbox{if } (m+k,r+s)\preccurlyeq (n,0)\\
   \langle\min\{(n,0),\big(3n-(m+k+l+1),-(r+s+t)\big)\},\bfm{0}\rangle, \hbox{ otherwise}
   \end{array}
 \right.\\
 && \mbox{since } (l,t)\preccurlyeq(n-1,0)\mbox{ and }(m+k,r+s)\preccurlyeq (n,0) \mbox{ imply}\\
 && \big(3n-(m+k+l+1),-(r+s+t)\big)\succcurlyeq (n,0), \mbox{ we have}\\
     &=&\langle \min\{(n,0),\big(3n-(m+k+l+1),-(r+s+t)\big)\},\bfm{0} \rangle,
 \eeqn
  so
$
\langle(m,r),\bfm{p}\rangle\odot\big[\langle(k,s), \beta\rangle\odot\langle(l,t),\gamma  \rangle\big]
 =\big[\langle(m,r),\bfm{p}\rangle\odot\langle(k,s), \beta\rangle\big]\odot\langle(l,t),\gamma  \rangle.
 $
  \subsubsection{$\bfm{0}\not=\beta\ast\gamma$}
 In this case  the associativity of $\odot$ also follows  from de associativity of $\ast$  in $\bfm{\LLn}$.

 \subsection{$\bfm{0}=\alpha$ and $\bfm{0}<_{\bfm{L_{p+1}}}\beta,\gamma$\label{A.3}}
 \subsubsection{$\alpha\ast\beta=\bfm{0}$}
 In this case  $\bfm{p}\notin\{\beta,\gamma\}$, and
 \beqn
  \lefteqn{\langle(m,r),\bfm{0}\rangle\odot\big[\langle(k,s), \beta\rangle\odot\langle(l,t),\gamma  \rangle\big]}\hspace{1cm}\\
  &=&\langle(m,r),\bfm{0}\rangle\odot\big\langle\min\{(n,0),(2n-(k+l+1),-(s+t)\},
 \bfm{0}  \rangle\\
   &=& \left\{\begin{array}{ll}
      \langle(m,r),\bfm{0}\rangle\odot \langle(n,0),\bfm{0}\rangle\, \ \hbox{if } (n,0)\preccurlyeq (2n-(k+l+1),-(s+t))\\
      \langle(m,r),\bfm{0}\rangle\odot \langle(2n-(k+l+1),-(s+t)),\bfm{0}\rangle, \quad  \hbox{otherwise}
       \end{array}
     \right.\\
 &=& \left\{\begin{array}{ll}
      \langle(n,0),\bfm{0}\rangle, \quad \ \hbox{if } (n,0)\preccurlyeq (2n-(k+l+1),-(s+t)))\\
      \langle\min\{(n,0),(a+[2n-(k+l+1)]+1,r-(s+t))\},\bfm{0} \rangle, \ \hbox{otherwise}
       \end{array}
     \right.\\
     && \mbox{since }(n,0)\preccurlyeq (2n-(k+l+1),-(s+t)) \mbox{ implies }\\
     &&(n,0)\preccurlyeq (m+[2n-(k+l+1)]+1,r-(s+t)),  \mbox{ we have}\\
 &=&  \min\{(n,0),\big(m+[2n-(k+l+1)]+1,r-(s+t)\big)\},\bfm{0} \rangle,\\
  &=&  \min\{(n,0),\big(m+2n-(k+l),r-(s+t)\big)\},\bfm{0} \rangle.
  \eeqn
  Moreover,
\beqn
\lefteqn{\big[\langle(m,r),\bfm{0}\rangle\odot\langle(k,s), \beta\rangle\big]\odot\langle(l,t),\gamma  \rangle}\hspace{1.25cm}\\
  &=&\langle(k,s)\to(m,r),\bfm{0}\rangle\odot\langle(l,t),\gamma \rangle\\
  &=&\langle\min\{(n,0),(n-k+m,-s+r),\bfm{0}\rangle\odot\langle(l,t),\gamma \rangle\\
  &=&\left\{\begin{array}{ll}
      \langle(n,0),\bfm{0}\rangle\odot \langle(l,t),\gamma\rangle\, \ \hbox{if } (n,0)\preccurlyeq (n-k+m,-s+r)\\
      \langle(n-k+m,-s+r),\bfm{0}\rangle\odot \langle(l,t),\gamma\rangle, \  \hbox{otherwise}
       \end{array}
     \right.\\
  &=&\left\{\begin{array}{ll}
      \langle(n,0),\bfm{0}\rangle, \ \hbox{if } (n,0)\preccurlyeq (n-k+m,-s+r)\mbox{ (Lemma~\ref{L: neu-ab})}\\
      \langle(l,t)\to(n-k+m,-s+r),\bfm{0}\rangle, \  \hbox{otherwise}
       \end{array}
     \right.\\
  &=&\left\{\begin{array}{ll}
      \langle(n,0),\bfm{0}\rangle, \ \hbox{if } (n,0)\preccurlyeq (n-k+m,-s+r)\ \\
      \langle\min\{(n,0),(n-l+n-k+m,-s+r-t)\},\bfm{0}\rangle, \  \hbox{otherwise}
       \end{array}
     \right.\\
  &=&\left\{\begin{array}{ll}
      \langle(n,0),\bfm{0}\rangle, \ \hbox{if } (n,0)\preccurlyeq (n-k+m,r-s)\preccurlyeq(0,0)\\
      \langle\min\{(n,0),(m+2n-(k+l),r-(s+t))\},\bfm{0}\rangle, \  \hbox{otherwise}
       \end{array}
     \right.\\
     && \mbox{since }(n,0)\preccurlyeq(n-k+m,r-s) \mbox{ and } (l,0)\preccurlyeq(n,0)\mbox{ imply }\\
     &&(n,0)\preccurlyeq (m+[2n-(k+l)],r-(s+t)),  \mbox{ we have}\\
  &=& \langle\min\{(n,0),(m+2n-(k+l),r-(s+t))\},\bfm{0}\rangle,
 \eeqn
 and so

$\langle(m,r),\bfm{0}\rangle\odot\big[\langle(k,s), \beta\rangle\odot\langle(l,t),\gamma  \rangle\big]=\big[\langle(m,r),\bfm{0}\rangle\odot\langle(k,s), \beta\rangle\big]\odot\langle(l,t),\gamma  \rangle.$\smallskip

 \subsubsection{$\alpha\ast\beta\not=\bfm{0}$}
  ${}$\vspace{-0.75cm}
  \bg{align*}
 \lefteqn{\langle(m,r),\bfm{0}\rangle\odot\big[\langle(k,s), \beta\rangle\odot\langle(l,t),\gamma  \rangle\big]}\hspace{1cm}\\ &=\langle(m,r),\bfm{0}\rangle\odot\big[\langle(k,s)\ast(l,t),
 \beta\ast\gamma  \rangle\big]\\
 &=\langle(m,r),\bfm{0}\rangle\odot\langle(k,s)\ast(l,t),\beta\ast\gamma\rangle
 = \langle([(k,s)\ast(l,t)]\to (m,r)), \bfm{0} \rangle\\
   &=\langle(k,s)\to((l,t)\to (m,r)), \bfm{0} \rangle
   =\langle(l,t)\to((k,s)\to (m,r)), \bfm{0} \rangle\\
   &=\langle(l,t),\gamma\rangle\odot\langle((k,s)\to (m,r), \bfm{0} \rangle
   =\langle(l,t),\gamma\rangle\odot[\langle(m,r), \bfm{0} \rangle\odot\langle(k,s),\beta\rangle]\\
   &=[\langle(m,r), \bfm{0} \rangle\odot\langle(k,s),\beta\rangle]\odot\langle(l,t),\gamma\rangle
 \end{align*}
\subsection{ $\bfm{0}<_{\bfm{L_{p+1}}}\alpha$ and $\beta=\gamma=\bfm{0}$.}
  ${}$\vspace{-0.75cm}
\beqn
 \lefteqn{\langle(m,r),\alpha\rangle\odot\big[\langle(k,s), \bfm{0}\rangle\odot\langle(l,t),\bfm{0}  \rangle\big]}\hspace{2cm}\\
 &=& \langle(m,r),\alpha\rangle\odot\langle\min\{(n,0),(k+l+1,s+t)\}, \bfm{0}\rangle\\
  &=&\left\{\begin{array}{ll}
      \langle(m,r),\alpha\rangle\odot \langle(n,0),\bfm{0}\rangle\, \ \hbox{if } (n,0)\preccurlyeq (k+l+1,s+t)\\
      \langle(m,r),\alpha\rangle\odot \langle(k+l+1,s+t),\bfm{0}\rangle, \quad  \hbox{otherwise}
       \end{array} \right.\\
  &=&\left\{\begin{array}{ll}
      \langle (n,0),\bfm{0}\rangle\, \ \hbox{if } (n,0)\preccurlyeq (k+l+1,s+t)\\
      \langle(m,r)\to(k+l+1,s+t),\bfm{0}\rangle, \quad  \hbox{otherwise}
       \end{array} \right.\\
       &=&\left\{\begin{array}{ll}
      \langle(n,0),\bfm{0}\rangle\, \ \hbox{if } (n,0)\preccurlyeq (k+l+1,s+t)\\
      \left\langle\min\left\{(n,0),\big((n-m)+(k+l+1),-r+(s+t)\big)\right\},\bfm{0}\right\rangle, \hbox{otherwise}
       \end{array} \right.\\
  && \mbox{since }(n,0)\preccurlyeq (k+l+1,s+t)\mbox{ and } (0,0)\preccurlyeq(n-m,-r) \mbox{ imply }\\
     &&(n,0)\preccurlyeq ((n-m)+k+l+1,-r+(s+t)),  \mbox{ we have}\\
  &=&\left\langle\min\left\{(n,0),\big((n-m)+k+l+1,-r+(s+t)\big)\right\}, \bfm{0}\right\rangle.
 \eeqn
 Moreover,
 \beqn
 \lefteqn{\big[\langle(m,r),\alpha\rangle\odot\langle(k,s), \bfm{0}\rangle\big]\odot\langle(l,t),\bfm{0}  \rangle}\hspace{2cm}\\
 &=& \langle(m,r)\to(k,s),\bfm{0}\rangle\odot\langle(l,t)\},\bfm{0}\rangle\\
 &=& \langle\min\{(n,0),(n-m+k,-r+s),\bfm{0}\rangle\odot\langle(l,t)\},\bfm{0}\rangle\\
  &=&\left\{\begin{array}{ll}
      \langle(n,0),\bfm{0}\rangle\odot \langle(l,t),\bfm{0}\rangle\, \ \hbox{if } (n,0)\preccurlyeq (n-m+k,-r+s)\\
      \langle(n-m+k,-r+s),\bfm{0}\rangle\odot \langle(l,t),\bfm{0}\rangle, \quad  \hbox{otherwise}
       \end{array} \right.\\
  &=&\left\{\begin{array}{ll}
      \langle(n,0),\bfm{0}\rangle \ \hbox{if } (n,0)\preccurlyeq (n-m+k,-r+s)\\
      \left\langle\min\{(n,0),\left((n-m+k)+l+1,(-r+s)+t\right)\},\bfm{0}\right\rangle, \quad  \hbox{otherwise}
       \end{array} \right.\\
  && \mbox{since }(n,0)\preccurlyeq (n-m+k,-r+t) \mbox{ implies }\\
     &&(n,0)\preccurlyeq ((n-m+k)+l+1,(-r+s)+t)  \mbox{ we have}\\
  &=& \left\langle\min\{(n,0),\left((n-m)+k+l+1,-r+(s+t)\right)\}. \bfm{0}\right\rangle.
 \eeqn
 And so
 \[\langle(m,r),\alpha\rangle\odot\big[\langle(k,s), \bfm{0}\rangle\odot\langle(l,t),\bfm{0}\rangle\big]=
 \big[\langle(m,r),\alpha\rangle\odot\langle(k,s), \bfm{0}\rangle\big]\odot\langle(l,t),\bfm{0}\rangle.\]
\subsection{$\alpha=\beta=\gamma=\bfm{0}$}
  ${}$\vspace{-0.75cm}
 \beqn
  \lefteqn{\langle(m,r),\bfm{0}\rangle\odot\big[\langle(k,s), \bfm{0}\rangle\odot\langle(l,t),\bfm{0}  \rangle\big]}\hspace{1.5cm}\\
 &=&\langle(m,r),\bfm{0}\rangle\odot\big\langle\min\{(n,0),(k+l+1,s+t)\},
 \bfm{0}  \rangle\\
 &=& \left\{\begin{array}{ll}
      \langle(m,r),\bfm{0}\rangle\odot \langle(n,0),\bfm{0}\rangle\, \ \hbox{if } (n,0)\preccurlyeq (k+l+1,s+t)\\
      \langle(m,r),\bfm{0}\rangle\odot \langle(k+l+1,s+t),\bfm{0}\rangle, \quad  \hbox{otherwise}
       \end{array}
     \right.\\
  &=& \left\{\begin{array}{ll}
       \langle(n,0),\bfm{0}\rangle\, \ \hbox{if } (n,0)\preccurlyeq (k+l+1,s+t)\\
       \langle\min\{(n,0),\big(m+(k+l+1)+1,r+(s+t)\big)\},\bfm{0}\rangle, \quad  \hbox{otherwise}
       \end{array}
     \right.\\
  && \mbox{since} (n,0)\preccurlyeq (k+l+1,s+t) \mbox{ implies }\\
  &&(n,0)\preccurlyeq (m+k+l+2,r+s+t), \mbox{ we have}\\
  &=&\langle\min\{(n,0),(m+k+l+2,r+s+t)\},\bfm{0}\rangle.
 \eeqn
 In similar way, and using the commutatiity of $\odot$ we have
 \begin{align*}
 \lefteqn{\big[\langle(m,r),\bfm{0}\rangle\odot\langle(k,s), \bfm{0}\rangle\big]\odot\langle(l,t),\bfm{0}\rangle}\hspace{2cm}\\
 & = \langle(l,t),\bfm{0}\rangle\odot\big[\langle(m,r),\bfm{0}\rangle
 \odot\langle(k,s), \bfm{0}\rangle\big]\\
 & =\langle\min\{(n,0),(m+k+l+2,r+s+t)\},\bfm{0}\rangle.
 \end{align*}
  ${}$\vspace{-0.75cm}

 And so
 \[\langle(m,r),\bfm{0}\rangle\odot\big[\langle(k,s), \bfm{0}\rangle\odot\langle(l,t),\bfm{0}\rangle\big] =\big[\langle(m,r),\bfm{0}\rangle\odot\langle(k,s), \bfm{0}\rangle\big]\odot\langle(l,t),\bfm{0}\rangle. \]

\subsection{$\bfm{0}<_{\bfm{L_{p+1}}}\alpha,\beta, \gamma$ }

 We consider subcases:
\subsubsection{$\alpha\ast(\beta\ast\gamma)\not=\bfm{0}$}

 This case is similar to \ref{A.1}  and \ref{A.2}, and the associativity of $\odot$ also follows  from de associativity of $\ast$  in $\bfm{\LLn}$.

 \subsubsection{$\alpha\ast(\beta\ast\gamma)=\bfm{0}$\label{A.6.2}}

 \subsubsection*{\ref{A.6.2}.1 $\alpha\ast\beta\not=\bfm{0}$ and $\beta\ast\gamma\not=\bfm{0}$}

 In this case $\alpha\not=\bfm{p}$. Then

 \beqn
 \lefteqn{\langle(m,r),\alpha\rangle\odot\big[\langle(k,s),\beta  \rangle\odot\langle(l,t),\gamma\rangle\big]}\hspace{0.25cm}\\
 &=&
  \langle(m,r),\alpha\rangle\odot\big\langle(k,s)\ast(l,t),\beta\ast\gamma\rangle,                                                                    \\&=&
  \langle(m,r),\alpha\rangle\odot\big\langle\max\{(0,0),(k+l-n,s+t),\beta\ast\gamma\rangle,                                                                    \\
 &=&\left\{\begin{array}{ll}
 \langle(m,r),\alpha\rangle\odot\langle(0,0),\beta\ast\gamma\rangle, \quad \hbox{ if } (k+l,t+s)\preccurlyeq(n,0), \\
  \langle\min\{(n,0),(2n-(m+(k+l-n)+1),-r-(s+t))\},\bfm{0}\rangle,  \hbox{ otherwise}
     \end{array}
   \right.\\
   &=&\left\{\begin{array}{ll}
 \langle\min\{(n,0),(2n-(m+1),r)\},\bfm{0}\rangle, \ \mbox{ if } (k+l,t+s)\preccurlyeq(n,0), \\
  \langle\min\{(n,0),(2n-(m+(k+l-n)+1),-r-(s+t))\},\bfm{0}\rangle,  \hbox{ otherwise}
     \end{array}
   \right.\\
   &=&\left\{\begin{array}{ll}
  \langle(n,0),\bfm{0}\rangle, \quad \hbox{ if } (k+l,t+s)\preccurlyeq(n,0), \\
  \langle\min\{(n,0),(3n-(m+k+l+1),-r-s-t)\},\bfm{0}\rangle,  \hbox{ otherwise}
     \end{array}
   \right.\\
    & &\mbox{ since }(m,r)\preccurlyeq(n-1,0) \mbox{ and }(k+l,t+s)\preccurlyeq(n,0)  \mbox{ implies }  \\&&(3n-(m+k+l+1),-r-s-t)\succcurlyeq(n,0), \mbox{ we have}\\
   &=& \langle\min\{(n,0),(3n-(m+k+l+1),-(r+s+t)\},\bfm{0}\rangle,
 \eeqn
 Moreover, by commutativity of $\odot$
 \beqn\lefteqn{\big[\langle(m,r),\alpha\rangle\odot\langle(k,s),\beta
 \rangle\big]\odot\langle(l,t),\gamma\rangle}\hspace{0.25cm}\\
 &=&
  \langle(l,t),\beta\rangle\odot \big[\langle(m,r),\alpha\rangle\odot\langle(k,s)\beta\rangle\big]                                          \\
   &=& \langle\min\{(n,0),(3n-(m+k+l+1),-(r+s+t)\},\bfm{0}\rangle.
  \eeqn
  And so
 \[
  \langle(m,r),\alpha\rangle\odot\big[\langle(k,s),\beta  \rangle\odot\langle(l,t),\gamma\rangle\big]=\big[\langle(m,r),\alpha\rangle\odot\langle(k,s),\beta  \rangle\big]\odot\langle(l,t),\gamma\rangle.
 \]
\subsubsection*{\ref{A.6.2}.2 $\alpha\ast\beta\not=\bfm{0}$ and $\beta\ast\gamma=\bfm{0}$}
 Similarly to the first equality in  \ref{A.2}.1, we  can see:
 \[
 \langle(m,r),\alpha\rangle\odot\big[\langle(k,s), \beta  \rangle\odot\langle(l,t),\gamma\rangle\big]
  = \langle\min\{(n,0),(3n-(m+k+l+1),-(r+s+t))\},\bfm{0}\rangle.
   \]
  Moreover, using the commutivity of $\odot$
 \beqn
  \lefteqn{\big[\langle(m,r),\alpha\rangle\odot\langle(k,s),\beta
 \rangle\big]\odot\langle(l,t),\gamma\rangle}\hspace{0.25cm}\\
  &=& \langle(l,t),\gamma\rangle\odot\langle(m,r)\ast(k,s),\alpha\ast\beta \rangle\\
  &=& \langle(l,t),\gamma\rangle\odot\langle\max\{(0,0),(m+k-n,r+s),\alpha\ast\beta \rangle\\
   &=&\left\{\begin{array}{ll}
 \langle(l,t),\gamma\rangle\odot\langle (0,0),\alpha\ast\beta\rangle, \ \mbox{ if } (m+k,r+s)\preccurlyeq(n,0), \\
  \langle(l,t),\gamma\rangle\odot\langle(m+k-n),r+s)\},\rangle, \alpha\ast\beta\rangle, \hbox{ otherwise}
     \end{array}
   \right.\\
   &=&\left\{\begin{array}{ll}
 \langle\min\{(n,0),(2n-(l+1),-t),\bfm{0}\rangle, \ \mbox{ if } (m+k,r+s)\preccurlyeq(n,0), \\
  \langle\min\{(n,0),(2n-(l+m+k-n+1),-(r+s)-t),\bfm{0}\rangle, \hbox{ otherwise}
     \end{array}
   \right.\\
   &=&\left\{\begin{array}{ll}
 \langle\min\{(n,0),(2n-(l+1),-t),\bfm{0}\rangle, \ \mbox{ if } (m+k,r+s)\preccurlyeq(n,0), \\
  \langle\min\{(n,0),(3n-(l+m+k+1),-(r+s+t),\bfm{0}\rangle, \hbox{ otherwise}
     \end{array}
   \right.\\
   &&\mbox{ since } (l,t)\preccurlyeq(n-1,0) \mbox{ imply } (n,0)\preccurlyeq(2n-(l+1),-t), \mbox{ we have,}\\
   &=&\langle\min\{(n,0),(3n-(l+m+k+1),-(r+s+t)),\bfm{0}\rangle.
   \eeqn
   And so
    \[
 \langle(m,r),\alpha\rangle\odot\big[\langle(k,s), \beta  \rangle\odot\langle(l,t),\gamma\rangle\big] =\big[\langle(m,r),\alpha\rangle\odot\langle(k,s),\beta
 \rangle\big]\odot\langle(l,t),\gamma\rangle.\]

 \subsubsection*{\ref{A.6.2}.3 $\alpha\ast\beta=\bfm{0}$ and $\beta\ast\gamma\not=\bfm{0}$}
 In this case the equality follows from  \ref{A.6.2}.2, using the commutativity of $\odot$.

 \bigskip

Any other case follows from the previous ones using the commutativity of $\odot$.

\section{Proof of Lemma \ref{L:*monoton3}}

 Let $a=\langle(m,r),\alpha\rangle$,  $b=\langle(k,s),\beta\rangle$ and $c=\langle(l,t),\gamma\rangle$ be arbitrary elements of $A^{p+1}_{n+1}$.
 We need to prove that  $b\le c$ implies $a\odot b\le a\odot c$.\smallskip

 Assume $b\le c$. Then $\beta\le_{L_{p+1}}\gamma$.
 \smallskip

 By definition of $\le$ and $\odot$,  the statement holds for $\alpha=\beta=\gamma=\bfm{p}$. We are going to prove it for the other cases.
 \be[\bf 1)]
 \item if \fbox{$\bfm{0}=\alpha=\beta=\gamma$}, then

  $b\le c$ implies $(m,r)\succcurlyeq(l,t)$, hence $(m+k+1,r+s)\succcurlyeq(m+l+1,r+t)$ so
  \bg{align*}
  a\odot b&= \langle\min\{(n,0),(m+k+1,r+s)\},\bfm{0}\rangle\\
  &\le \langle\min\{(n,0),(m+l+1,r+s)\},\bfm{0}\rangle=a\odot c.
  \end{align*}
  \item If \fbox{$\bfm{0}=\alpha$ and $\bfm{0}<_{\bfm{L_{p+1}}}\beta,\gamma$}, then

   $b\le c$ implies $(m,r)\le_{\LLn} (l,t)$, hence $(k,s)\to(m,r)\le_{\LLn} (l,t)\to(m,r)$ and so
  \[
  a\odot b= \langle(k,s)\to(m,r),\bfm{0}\rangle
  \le \langle(l,t)\to(m,r),\bfm{0}\rangle =a\odot c.
  \]
  \item If \fbox{$\bfm{0}=\alpha=\beta$ and $\bfm{0}<_{\bfm{L_{p+1}}}\gamma$}, then

  $b\le c$ implies $(n-1,0)\le_{\LLn}(k+l,s+t)$,
  hence $(n-l,-t)\le_{\LLn}(m+1,s)$ and so
 \bg{align*}
  a\odot b&= \langle\min\{(n,0),(m+k+1,r+s)\},\bfm{0}\rangle\\
  &\le \langle\min\{(n,0),(n-l+m,r-t)\},\bfm{0}\rangle\\
  &\le \langle(l,t)\to(m,r),\bfm{0}\rangle=a\odot c.
  \end{align*}
 \item  If \fbox{$\bfm{0}\notin\{\alpha,\beta,\gamma\}$}, then

 $b\le c$ implies  $(k,s)\le_{\LLn} (l,t)$.

 If \underline{$\alpha\ast\beta\not=\bfm{0}$}, then $\alpha\ast\gamma\not=\bfm{0}$, hence

 $(m,r)\ast (k,r)\le_{\bfm{\LLn}}(m,r)\ast (l,t)$
  and so

 $a\odot b=\langle(m,r)\ast (k,r),\alpha\ast\beta\rangle\le \langle(m,r)\ast (l,t),\alpha\ast\gamma\rangle= a\odot c.$

 If \underline{$\alpha\ast\gamma=\bfm{0}$}, then $\alpha\ast\beta=\bfm{0}$. hence\\
 $(2n-(m+k+1),-(r+s))\succcurlyeq(2n-(m+l+1),-(r+t))$,
  and so
  \bg{align*}
  \lefteqn{\min\{(n,0),(2n-(m+k+1),-(r+s))\}}\hspace{1.5cm}\\
  &\le_{\LLn}\min\{(n,0),(2n-(m+l+1),-(r+t))\}.
  \end{align*}
 Thus $a\odot b\le a\odot c.$

 If \underline{$\alpha\ast\beta=\bfm{0}$} and \underline{$\alpha\ast\gamma\not=\bfm{0}$}
  \begin{align*}
      a\odot b & = \langle \min\{(n,0),(2n-(m+k+1),-(r+s))\},\bfm{0}\rangle,\\   a\odot c &=\langle \max\{(0,0),(m+l-n,r+t)\},\alpha\ast\gamma\rangle
 \end{align*}
 If $a\odot b=\langle (n,0),\bfm{0}\rangle$, then $a\odot b\le a\odot c$.\\
 For $a\odot b= \langle (2n-(m+k+1),-(r+s)),\bfm{0}\rangle$, we have two possibilities
 \bi
  \item $a\odot c =\langle (m+l-n,r+t),\alpha\ast\gamma\rangle $ and so
 \begin{align*}
  \lefteqn{([2n-(m+k+1)]+[m+l-n],-(r+s)+(r+t))}\hspace{3cm}\\
  &=(n+(l-m)+1),-s+t)\\
  &= (n-1,0)+(l-m,t-s)\\ &\succcurlyeq (n-1,0)+(0,0)
  \succcurlyeq (n-1,0)
 \end{align*}
  and  $a\odot b\le a\odot c.$
\item $a\odot c =\langle (0,0),\alpha\ast\gamma\rangle $ In this case  $(m+k,r+s)\preccurlyeq (m+l,r+t)\preccurlyeq (n,0)$. Hence
    \begin{align*}
 \lefteqn{(2n-(m+k+1),-(r+s))}\hspace{3cm}\\
  &=(2n-1,0)-(m+k,r+s)\succcurlyeq\\
  &= (2n-1,0)-(n,0)= (n-1,0),
 \end{align*}
 and so $a\odot b\le a\odot c.$
\ei

   \ee
The other cases can be deduced from the previous ones \hfill$\Box$

\end{document}